\documentclass[10pt]{article}
\title{Uniform bounds for point cohomology of $\lp{1}(\Z_+)$ and related algebras}

 \author{\sc Yemon Choi\\
}
 \date{\sl\small Dedicated to the memory of Graham R.~Allan, \textup{1936--2007}
}

\RequirePackage{amsmath,amssymb,amsfonts,amsbsy}

\usepackage[usenames]{color}

\RequirePackage[PS,nohug]{diagrams}
\newcommand{\RL}[2]{\pile{\rTo^{#1} \\ \lTo_{#2} }} 

\newcounter{pulse}[section]
\numberwithin{pulse}{section}  

\newcounter{alph}

\numberwithin{equation}{section}

\newcommand{\Cplx}{\mathbb C}
\newcommand{\Nat}{\mathbb N}
\newcommand{\Rat}{\mathbb Q}
\newcommand{\Real}{\mathbb R}
\newcommand{\Z}{{\mathbb Z}}

\newcommand{\Disc}{\mathbb D} 
\newcommand{\bbG}{\mathbb G} 
\newcommand{\UHP}{{\mathbb H}} 
\newcommand{\Torus}{\mathbb T} 

\newcommand{\AT}{A({\mathbb T})}  
\newcommand{\ATpos}{A({\mathbb T})^+}  

\newcommand{\st}{\;:\;}
\newcommand{\defeq}{:=}
\newcommand{\sid}{\mathop{\sf id}}  

\newcommand{\veps}{\varepsilon}
\newcommand{\wtild}{\widetilde}

\newcommand{\blank}{\underline{\quad}}
\newcommand{\iso}{\cong}

\renewcommand{\Im}{\mathop{\rm Im}}

\newcommand{\id}[1][]{{\sf 1}_{#1}} 
\newcommand{\wstar}{\ensuremath{{\rm w}^*}}

\newcommand{\tp}{\mathop{\otimes}\nolimits}
\newcommand{\ptp}{\mathop{\widehat{\otimes}}\nolimits} 

\newcommand{\abs}[1]{\vert{#1}\vert}
\newcommand{\norm}[1]{\Vert{#1}\Vert}
\newcommand{\Norm}[1]{\left\Vert{#1}\right\Vert}

\newcommand{\clos}[1]{\overline{#1}}

\newcommand{\Lin}[2]{{\mathcal L}({#1},{#2})}
\newcommand{\lp}[1]{\ell^{#1}} 

\newcommand{\dif}{\delta}
\newcommand{\Co}[3][]{\mathcal{#2}^{#3}{#1}}  
\DeclareMathOperator{\Ext}{Ext}

\newcommand{\cD}{{\mathcal D}}
\newcommand{\cF}{{\mathcal F}}
\newcommand{\bF}{{\mathbf F}}
\newcommand{\cG}{{\mathcal G}}
\newcommand{\bG}{{\mathbf G}}
\newcommand{\bj}{\mathop{\mathbf j}\nolimits}
\newcommand{\cJ}{{\mathcal J}}
\newcommand{\cL}{{\mathcal L}}
\newcommand{\bM}{{\mathbf M}}
\newcommand{\cM}{{\mathcal M}}

\newcommand{\ssp}{{\sf p}}
\newcommand{\ssS}{{\sf S}}
\newcommand{\bT}{{\mathbf T}}

\newcommand{\al}{\alpha}
\newcommand{\gm}{\gamma}
\newcommand{\lm}{\lambda}
\newcommand{\lmbar}{\bar{\lm}}   
\newcommand{\Lm}{\Lambda}

\newcommand{\indic}[1]{{\bf 1}_{[#1]}}

\newcommand{\cu}[1]{{#1}_{\rm un}}  
\newcommand{\fu}[1]{{#1}^{\sharp}}

\newcommand{\dual}[1]{\widehat{#1}}

\newcommand{\upos}[1]{{u_{#1}^\oplus}}
\newcommand{\uneg}[1]{{u_{#1}^\ominus}}

\newcommand{\dt}[1]{\textcolor{Bittersweet}{\sf#1}}

\renewcommand{\emph}[1]{{\sl #1\/}} 

\newenvironment{acknowledgements}
{\subsection*{Acknowledgements}}
{\ignorespacesafterend}

\RequirePackage{amsthm}   

\theoremstyle{plain}
\newtheorem{thm}[pulse]{\sc Theorem}
\newtheorem{propn}[pulse]{\sc Proposition}
\newtheorem{lem}[pulse]{\sc Lemma}
\newtheorem{coroll}[pulse]{\sc Corollary}
\newtheorem{qu}[alph]{Question}
\theoremstyle{definition}
\newtheorem{defn}[pulse]{\sc Definition}
\theoremstyle{remark}
\newtheorem{rem}[pulse]{\sc Remark}

\newenvironment{YCqu}
{\begin{qu}}
{\end{qu}\ignorespacesafterend}

\newenvironment{YCproof}
{\begin{proof}}
{\end{proof}\ignorespacesafterend}

\newenvironment{YCnameproof}[1]
{\begin{proof}[{#1}]}
{\end{proof}\ignorespacesafterend}

\begin{document}
\maketitle

\begin{abstract}
It is well-known that the point cohomology of the convolution algebra $\ell^1({\mathbb Z}_+)$ vanishes in degrees 2 and above. We sharpen this result by obtaining splitting maps whose norms are bounded independently of the choice of point module. Our construction is a by-product of new estimates on projectivity constants of maximal ideals in $\lp{1}(\Z_+)$\/. Analogous results are obtained for some other $L^1$-algebras which arise from `rank one' subsemigroups of ${\mathbb R}_+$\/.
\end{abstract}



\begin{section}{Introduction}
If $A$ is a commutative, semisimple Banach algebra, the simplest $A$-modules are the 1-dimen\-sional modules corresponding to the characters of $A$\/, sometimes called \dt{point modules} as they correspond to evaluation of Gelfand transforms of elements of~$A$.
The calculation of cohomology with coefficients in such modules (henceforth referred to as \dt{point cohomology}) has been studied by several authors, with particular emphasis on the case where $A$ is a uniform algebra (see \cite{Erm_H2,Hel_UAptcoho}). 
Part of the motivation comes from classical commutative algebra, where point derivations on the coordinate ring ${\sf R}$ of an affine variety ${\mathcal V}\subset k^n$ encode \emph{intrinsic} geometric information about ${\mathcal V}$\/.
(More precisely: if $x\in{\mathcal V}$, then the space of $k$-linear derivations ${\sf R} \to k_x$ may be identified with the \dt{tangent space} of $V$ at $x$.)

In this article we instead focus on certain $\lp{1}$-analogues of the disc algebra.
The starting point for our investigations was the following question, which the author learned of from M.~C. White.
Consider the convolution algebra $\lp{1}(\Z_+)$, whose maximal ideal space is naturally identified with the closed unit disc $\clos{\Disc}$; for each $\lm\in\clos{\Disc}$ let $\Cplx_\lm$ denote the corresponding point module.

\begin{YCqu}\label{mainQ}
Do there exist constants $(C_n)_{n\geq 1}$ such that, for any $\lm\in\clos{\Disc}$, the coboundary map $\delta:\Co{C}{n}(\lp{1}(\Z_+),\Cplx_\lm)\to\Co{C}{n+1}(\lp{1}(\Z_+),\Cplx_\lm)$ is open with constant $< C_n$\/?
\end{YCqu}

If we allow the constants $C_n$ to vary with $\lm$ then the answer is known to be yes: the point is that we want to solve cohomology problems with norm control that is \emph{uniform in the choice of point module}.
%
Even the case $n=1$ appears to have gone unsolved thus far.

\begin{rem}
In the case $n=1$, Question~\ref{mainQ} is as follows: does there exist $C_1>0$ such that, whenever $\psi\in \lp{\infty}(\Z_+)$ and $\lm\in\clos{\Disc}$ satisfy
\[   \abs{\lm^j\psi_k -\psi_{j+k} + \psi_j \lm^k} \leq 1 \qquad\text{ for all $j,k\in\Z_+$}\;, \] 
there exists $\alpha\in\Cplx$ such that
\[  \abs{\psi_n - n\lm^{n-1}\alpha} \leq C_1 \qquad \text{ for all $n\in\Nat$\/?} \]
More informally: is every approximate point derivation on $\lp{1}(\Z_+)$ \emph{uniformly} close to a true point derivation? 
If we were to allow $C_1$ to vary with $\lm$ then, by a trivial inductive argument, taking $C_1=\sup_n n\abs{\lm}^{n-1} = O(1-\abs{\lm})^{-1})$ would suffice.
\end{rem}

In this article we will answer Question~\ref{mainQ} in the affirmative (Corollary~\ref{c:stable-cocycle} below), by giving explicit splitting maps for the Hochschild chain complex; moreover, the constants $C_n$ can all be
bounded independently of~$n$\/, though we do not know if our estimates are sharp.
Our approach is influenced by arguments of Johnson in which he calculates the point cohomology of the disc algebra (see \cite[Proposition 9.1]{BEJ_CIBA}); but we also draw on ideas of \dt{projective} and \dt{flat} Banach modules, in the sense of Helemskii. In fact, the desired uniform estimates for our splitting maps follow from an analysis of the \dt{projectivity constants} of the maximal ideals in $\lp{1}(\Z_+)$\/, and a~secondary purpose of this article is to lend a quantitative flavour to some otherwise-familiar arguments from the homo\-logical theory of Banach modules (see also \cite[\S\S1--4]{MCW_UA}).
The methods we use extend to the convolution algebras of some other totally ordered, cancellative  abelian semigroups. We can obtain similar results for $\lp{1}(\bbG_+)$\/, where $\bbG$ is a dense subgroup of $\Rat$\/, and for $L^1(\Real_+)$\/: higher rank examples are beyond the scope of the present work, and we hope to address them in a sequel.

\begin{rem}
Uniform bounds are easily obtained -- and are well-known -- if we restrict attention to characters which correspond to \dt{peak points} of $\lp{1}(\Z_+)$\/, i.e.~those corresponding to the unit circle.
The force of our result is to produce uniform bounds for the point modules which arise from~$\Disc$\/.
Such uniform bounds are to be expected for the disc algebra, whose group of isometric automorphisms acts \emph{transitively} on $\Disc$,
and which therefore ought not to privilege any
point of $\Disc$ over any other. In contrast, an elementary argument using extreme points shows that
every isometric automorphism of $\lp{1}(\Z_+)$ is given by rotation of $\Disc$\/: so the achievement of uniform bounds over \emph{all} point modules does not follow from trivial symmetry arguments.
\end{rem}

\begin{rem}
Interestingly, although there is {\it a priori} no relation between the point cohomology of a commutative  Banach algebra and that of its uniform completion, the splitting maps we construct will also work just as well for the uniform completions of $\lp{1}(\Z_+)$\/, $\lp{1}(\bbG_+)$ and $L^1(\Real_+)$\/. The significance of this phenomenon is not clear to the author at present, and may warrant further attention.
\end{rem}
\end{section}

\begin{section}{Notation and preliminaries}
We start by fixing some notation. Throughout we shall denote the identity map on a Banach space, module or algebra by $\sid$ (it will be clear from context what the domain of $\sid$~is). The dual of a Banach space $E$ will be denoted throughout by $E'$\/, and the adjoint of a bounded linear map $f$ by~$f'$\/.
The canonical (isometric) embedding of $E$ into its second dual $E''$ will be denoted by~$\kappa_E$\/, and $\ptp$ will denote the projective tensor product of Banach spaces.

We refer to standard references such as \cite{Bons-Dunc} for the definitions of Banach algebras and bimodules over them.
(In particular, bimodule actions are \emph{not} assumed {\it a~priori}\/ to be contractive.)
The \dt{forced unitization} $A\oplus\Cplx$ of a Banach algebra (see \cite[Definition 3.1]{Bons-Dunc}) will be denoted by $\fu{A}$\/: note that we are taking the $\lp{1}$-sum of $A$ and $\Cplx$\/.

\begin{rem}
In several places below we will consider bounded nets in $\Lin{V}{W'}$\/, where $V$ and $W$ are Banach spaces, and speak of taking weak$^*$-cluster points in $\Lin{V}{W'}$\/. The topology referred to is the weak$^*$-topology conferred on $\Lin{V}{W'}$ by its canonical isometric predual $V\ptp W$\/, and \emph{not} the weak$^*$-operator topology defined by the functionals $\theta_{v,w}:T\mapsto T(v)(w)$\/.
\end{rem}

\begin{subsection}{Certain convolution algebras}\label{s:L1prelim}
The examples we consider are all of the form $L^1(G_+)$\/, where $G$ is a locally compact group that admits a continuous injective homomorphism $\iota:G\to\Real$\/; $G_+$ is defined to be the inverse image under $\iota$ of the non-negative reals. The study of such algebras goes back to Arens and Singer~\cite{AreSin_gaf56}, who also studied the completions of these algebras in the uniform norm to obtain analogues of the disc algebra.

We start with some general preliminary results. First some notation:  $\dual{G}$ denotes the set of all \dt{characters} of $G$\/, that is, all the continuous group homomorphisms from $G$ into~$\Torus$\/; and $\mu$ will be a fixed Haar measure on~$G$\/.

\begin{lem}\label{l:reiter}
There exists a net $(u_\al)$ of positive, compactly supported, integrable functions on $G_+$\/, such that $\int_G u_\al \,d\mu=1$ for all $\al$, and
\[ \norm{u_\al*f-u_\al\cdot \bigl( \int_G f\,d\mu \bigr) } \to 0 \quad\text{ for all $f\in L^1(G_+)$\/.} \]
\end{lem}

\begin{YCproof}
Since every abelian group satisfies Reiter's $(P_1)$-condition,
there exists a net $(v_\al) \subset L^1(G)$ with $\nu_\al\geq 0$, $\int_G \nu_\al\,d\mu=1$ and $\norm{\nu_\al*\delta_x-\nu_\al}\to 0$ uniformly on compact subsets of~$G$\/. The construction is easily modified to ensure that each $\nu_\al$ has compact support, and a standard approximation argument shows that
\[ \norm{\nu_\al*f-\nu_\al \cdot\bigl( \int_G f\,d\mu \bigr) } \to 0 \quad\text{ for all $f\in L^1(G_+)$\/.} \]
For each $\al$\/, since each $\nu_\al$ has compact support and $\iota(G)\subset\Real$, there exists $y(\al)\in G_+$ such that $y(\alpha)+\mathop{\rm supp}(\nu_\al)\subset G_+$\/. Putting $u_\al=\delta_{y(\al)}*\nu_{\al}$\/, we find that the net $(u_\al)$ has the required properties.
\end{YCproof}

\begin{coroll}\label{c:making-delta-nets}
Let $\gm\in\dual{G}$\/. Then there exists a net $(t_\al)\subset L^1(G_+)$ such that $\norm{t_\al}=1$ and $\norm{t_\al*f -t_\al\cdot\hat{f}(\gm)} \to 0 $\/.
\end{coroll}

\begin{YCproof}
We may exploit the action of the dual group $\dual{G}$ on $L^1(G)$\/, to transfer results for the augmentation character to the character~$\gm$\/.
Explicitly: let $(u_\al)$ be the net provided by Lemma~\ref{l:reiter} and put $t_\al(x)=u_\al(x)\gm(x)^{-1}$, $x\in G$\/. Given $f\in L^1(G_+)$\/, let $h(x)=f(x)\gm(x)$ for all $x\in G$. Then $h\in L^1(G)$\/, and
$(t_\al*f)(y) = (u_\al*h)(y)\gm(y)^{-1}$
for all $y\in G$\/;
also,
\[ \hat{f}(\gm) = \int_G f(x)\gm(x)\,d\mu(x) = \int_G h \,d\mu\,.\]
It follows from Lemma~\ref{l:reiter} that
\[ \norm{t_\al*f-t_\al\cdot\hat{f}(\gm)} = \norm{u_\al*h - u_\al\cdot \bigl(\int_G h\,d\mu\bigr)} \to 0 \;. \]
\end{YCproof}

In this article we restrict attention to the following three cases: $G=\Z$\/; $G$ a subgroup $\bbG\subseteq \Rat$\/; and $G=\Real$ with its usual topology.
In each case, to calculate the point cohomology we must first determine the set of characters. The character space of $\lp{1}(\bbG_+)$ will be addressed in Section~\ref{ss:l1Q+}.
The characters on $\lp{1}(\Z_+)$ and $L^1(\Real_+)$ are well-known, so we shall 
merely fix our notation: details can be found for instance, in \cite[Theorems 4.6.9 and 4.7.27]{Dal_BAAC}. We have already used $\Disc$ to denote the open unit disc in~$\Cplx$\/; the open upper half-plane $\{z\st \Im z >0\}$ will be denoted by~$\UHP$\/, and its closure in $\Cplx$ by~$\clos{\UHP}$\/.

The character space of $\lp{1}(\Z_+)$ may be identified with the closed unit disc $\clos{\Disc}$\/, and the Gelfand transform $\cG: \lp{1}(\Z_+)\to C(\clos{\Disc})$ is given by
\begin{equation}
(\cG a)(z) = \sum_{n\in\Z_+} a_n z^n \qquad(z\in \clos{\Disc})\,.
\end{equation}
Thus, the range of $\cG$ is the algebra $\ATpos$ consisting of all analytic functions on $\Disc$ whose Taylor series (about~$0$) is absolutely convergent.
Alternatively, we may identify $\ATpos$ with the space of all $f\in\AT$ such that
$\widehat{f}(n)=0$ for all $n < 0$\/,
where as usual $\AT$ denotes the algebra of functions on $\Torus$ with absolutely convergent Fourier series.

Let $\cF: L^1(\Real)\to C_0(\Real)$ denote the \dt{Fourier transform}
\begin{equation}
 (\cF f)(x) = \int_{-\infty}^\infty f(t)e^{itx}\,dt \,.
\end{equation}
If $f\in L^1(\Real_+)$ then $\cF f$ extends continuously to $\clos{\UHP}$  and is holomorphic on $\UHP$\/: the extension
$\cL f$ is essentially given by the \dt{Laplace transform} of~$f$\/. Every character on $L^1(\Real_+)$ is of the form $\cL(\blank)(\lm)$ for some $\lm \in\clos{\UHP}$\/; if $\lm\in\Real$\/, then the corresponding character extends to a character on $L^1(\Real)$\/.
\end{subsection}

\begin{subsection}{Cohomology, projectivity and flatness}
The definitions of Hochschild cohomology groups for Banach algebras and modules
can be found in the standard sources \cite{Hel_HBTA,BEJ_CIBA}\/:
we follow the notation of
of \cite[\S I.3.1]{Hel_HBTA}. Thus, given a Banach algebra $A$ and a Banach $A$-bimodule~$M$\/: the space of (continuous) \dt{$n$-cochains} is denoted by $\Co{C}{n}(A,M)$\/; the \dt{Hochschild coboundary operator} by $\dif$\/; and the (continuous) \dt{Hochschild cohomology groups} by $\Co{H}{n}(A,M)$\/, $n=0,1,2,\ldots$
%


Let $M_\lm =\{ a\in \ell^1(\Z_+)\;:\; \sum_{n\geq 0} a_n\lm^n=0\}$\/; this is a maximal ideal in $\ell^1(\Z_+)$.
 The existing proofs that $\Co{H}{n}(\lp{1}(\Z_+),\Cplx_\lm)=0$ for all $n\geq 2$ and all $\lm\in\clos{\Disc}$ are based, implicitly or explicitly, on the following two observations:
\begin{enumerate}
\item if $\abs{\lm}=1$ then $\Cplx_\lm$ is \dt{flat} as a left (or right) Banach $\lp{1}(\Z_+)$-module\/;
\item if $\abs{\lm}<1$ then $M_\lm$ is \dt{projective}
 as a left (or right) Banach $\lp{1}(\Z_+)$-module.
\end{enumerate}
This is part of more general theory, in which one uses homological properties such as \dt{projectivity} and \dt{flatness} to prove vanishing theorems for Hochschild cohomology, often with the implicit construction of splitting maps. For Banach modules there are natural \emph{quantitative} variants, which can be used to control of the norms of these splitting maps.

The first systematic account appears to be in~\cite{MCW_UA}: we briefly review the material we will need.
It is convenient to introduce the \dt{conditional unitization} of a Banach algebra $A$\/: this is defined to be $A$ itself if $A$ is unital, and $\fu{A}$ otherwise.
Every left, right or two-sided Banach $A$-module is naturally a Banach $\cu{A}$-module of the respective type.

The following definitions are not the original ones but are equivalent to them. (We work with right modules rather than left modules as this will be marginally more convenient for us in Propositions \ref{p:proj-gen-split} and \ref{p:flat-gen-split} below.)
Throughout this section and what follows, we use $\pi$ to denote the continuous linear map $M\ptp\cu{A} \to M$ that is given by $\pi(x\tp a)=xa$ for all $x\in M$ and $a\in \cu{A}$\/.

\begin{defn}[cf.~{\cite[Proposition~2.8]{MCW_UA}}]
\label{dfn:proj-with-constant}
Let $A$ be a Banach algebra and $M$ a right Banach $A$-module, and let~$C>0$\/. We say that $M$ is \dt{$A$-projective with constant $\leq C$} if there exists a bounded linear, right $A$-module map $\sigma: M\to M\ptp\cu{A}$  such that $\pi\sigma=\sid$ and $\norm{\sigma}\leq C$\/.
\end{defn}

\begin{defn}[cf.~{\cite[Propositions 3.8 and 4.9]{MCW_UA}}]
\label{dfn:flat}
Let $A$ be a  Banach algebra and $M$ a right Banach $A$-module, and let~$C>0$\/. We say that $M$ is \dt{$A$-flat with constant $\leq C$} if there exists a bounded linear, left $A$-module map $\Lambda:  (M\ptp\cu{A})'\to M'$  such that $\Lambda\pi'=\sid$ and $\norm{\Lambda}\leq C$\/.
\end{defn}

In particular, an $A$-module which is projective with constant $\leq C$ is {\it a fortiori}\/ flat with constant $\leq C$\/ (just take $\Lm=\sigma'$).


We finish this section with a lemma that, while not needed for the solution of Question~\ref{mainQ}, will be useful when considering the point cohomology of more general examples. It gives a sufficient criterion for flatness that is surely not new, but for which we have been unable to find a precise reference in the literature.

\begin{lem}\label{l:flat-by-approx}
Let $A$ be a Banach algebra and let $M$ be a right Banach $A$-module. Suppose there exists a bounded net $(\rho_\al)$ of bounded linear maps $M\to M\ptp \cu{A}$ such that
\begin{enumerate}
\item for every $x\in M$\/, $\pi\rho_\al(x)$ converges weakly to~$x$\/;
\item for every $x\in M$ and $a\in \cu{A}$\/, $\rho_\al(x\cdot a)-\rho_\al(x)\cdot a$ converges weakly to~$0$\/.
\end{enumerate}
Then $M$ is $A$-flat with constant $\leq \sup_\al \norm{\rho_\al}$\/.
\end{lem}

\begin{YCproof}%
Take $\Lm$ to be a \wstar-cluster point of the net $({\rho_\al}')\subset\Lin{(M\ptp\cu{A})'}{M'}$\/.
Clearly $\norm{\Lambda}\leq\sup_\al\norm{\rho_\al}$\/.
Routine calculations show that $\Lambda\pi'=\sid$ (using the first condition) and that $\Lambda(a\cdot\Psi)=a\cdot\Lambda(\Psi)$ for all
 $\Psi\in (M\ptp\cu{A})'$ and $a\in A$ (using the second condition).
%
\end{YCproof}

\end{subsection}
\end{section}

\begin{section}{Constructing splitting maps with controlled norms}
Throughout this section $A$ denotes a fixed commutative semisimple Banach algebra, not necessarily with identity. Let $\varphi$ be a character on~$A$ and let $M_\varphi=\ker\varphi$ be the corresponding maximal ideal. We shall construct splitting maps for (portions of) the Hochschild chain complex $\Co{C}{*}(A,\Cplx_\varphi)$\/, and estimate their norms, under various hypotheses on $\Cplx_\varphi$ and $M_\varphi$\/.

In the main applications below, $A$ has a bounded approximate identity and $M_\varphi$ is $A$-flat.
If we were to use the machinery of $\Ext$ for Banach modules, then we could easily deduce from these assumptions that certain cohomology groups are zero,
as follows:
for each $n\geq 1$ there are isomorphisms
\[ \Co{H}{n}(A,\Cplx_\varphi)\iso \Ext^n_A(\Cplx_\varphi,\Cplx_\varphi) = \Ext^{n-1}_A(M_\varphi,\Cplx_\varphi)\,, \]
and flatness of $M_\varphi$ implies that the right-hand side vanishes for $n\geq 2$\/.

This argument does not in itself provide the splitting maps we require. Nevertheless, the map $\Lambda$ which demonstrates flatness of $M_\varphi$
can be chased through the homological machinery of `comparison of resolutions', until one arrives at explicit splitting maps for the Hochschild chain complex $\Co{C}{*}(A,\Cplx_\varphi)$ in degrees~$\geq 2$\/.
Rather than giving the details in full generality, we shall merely
write down appropriate splitting maps -- as if by {\it ad hoc}\/ construction -- rather than showing how one is led to them. This keeps our account shorter,
albeit at the expense of some motivation for various formulas.

We start with the projective case, since \emph{this will be sufficient to resolve Question~\ref{mainQ}}, and since it provides a guiding outline for the proof of the flat case.
Before doing so, it is useful to make the following definition: if $T\in\Co{C}{n}(A,M)$ then $T$ is a bounded multi\-linear map from the $n$-fold Cartesian product $A\times\ldots \times A$ to~$M$\/; so by the universal property of the projective tensor product, we may identify $T$ with a unique
bounded linear map from $A^{\ptp n}$ to $M$ (the identification preserves norms). We shall denote this linear map by~$\wtild{T}$\/.

\begin{propn}\label{p:proj-gen-split}
Let $A$ be a unital commutative Banach algebra and let $\varphi\in\Phi_A$\/. Let $X$ be a Banach $A$-bimodule where the left action of $A$ is given by
\[  f\cdot x = \varphi(f)x \qquad(f\in A, x\in X). \]
Suppose that $M_\varphi$ is $A$-projective with constant~$\leq C$\/.
Then there exist bounded linear maps $s_1:\Co{C}{2}(A,X)\to \Co{C}{1}(A,X)$ and $s_2:\Co{C}{3}(A,X)\to\Co{C}{2}(A,X)$\/, such that $s_2\dif+\dif s_1=\sid$\/.
\[ 
\begin{diagram}[tight,width=4em]
\Co{C}{1}(A,X) & \RL{\dif}{s_1} & \Co{C}{2}(A,X) & \RL{\dif}{s_2} & \Co{C}{2}(A,X)
\end{diagram}
\]
Moreover, we can choose $s_1$ and $s_2$ such that $\max(\norm{s_1}, \norm{s_2}) \leq 2C+1$\/.
\end{propn}

\begin{YCproof}
By hypothesis, there exists a bounded linear, right $A$-module map $\sigma: M_\varphi\to M_\varphi\ptp A$ such that $\norm{\sigma}\leq C$ and $\pi\sigma=\sid$\/. Define $\ssp: A\to M_\varphi$ by $\ssp(f)=f-\varphi(f)\id$\/, and define bounded linear maps $s_1:\Co{C}{2}(A,X)\to\Co{C}{1}(A,X)$ and $s_2:\Co{C}{3}(A,X)\to\Co{C}{2}(A,X)$ by
\begin{flalign*}
\qquad (s_1 F)(f) & = -\wtild{F}\sigma\ssp(f) + \varphi(f)F(\id,\id) & (F\in\Co{C}{2}(A,X); f\in A);\qquad  \\
\qquad (s_2G)(f,g) & = -\wtild{F}(\sigma\ssp(f)\tp g) + \varphi(f)G(\id,\id,g) & (G\in\Co{C}{3}(A,X); f,g\in A). \qquad
\end{flalign*}

Note that if $w\in M_\varphi\ptp A$\/ and $b\in A$\/, we have
\begin{equation}\label{eq:GOOFY}
 \wtild{\dif F}(w\tp b) = - \wtild{F}(\pi(w)\tp b) + \wtild{F}(w\cdot b) -\wtild{F}(w)b
\end{equation}
(this is easily verified by approximating $w$ with a finite sum of elementary tensors in $M_\varphi\tp A$).
Let $f,g\in A$\/: then since $\varphi(fg)=\varphi(f)\varphi(g)$, we have
\begin{equation}\label{eq:BODGER}
\begin{aligned}
\dif s_1 F(f,g) & = \varphi(f)s_1 F(g) - s_1 F(fg) + s_1 F(f)g \\
	& = - \varphi(f) \wtild{F}\sigma\ssp(g) + \wtild{F}\sigma\ssp(fg) - \wtild{F}\sigma\ssp(f)g + \varphi(f) F(\id,\id)g \;;
\end{aligned}
\end{equation}
and, since $\sigma\ssp(f) \in M_\varphi\ptp A$\/, applying the identity \eqref{eq:GOOFY} gives
\begin{equation}\label{eq:BADGER} 
\begin{aligned}
s_2\dif F(f,g) & = -\wtild{\dif F}(\sigma\ssp(f)\tp g) + \varphi(f)\dif F(\id,\id, g) \\
	& = 	\wtild{F}(\ssp(f)\tp g) - \wtild{F}(\sigma\ssp(f)g)
		+ \wtild{F}\sigma\ssp(f)g
		+ \varphi(f) F(\id,g) - \varphi(f) F(\id,\id)g
\\
	& = F(f,g) - \wtild{F}\sigma(\ssp(f)g)  
		 + \wtild{F}\sigma\ssp(f)g- \varphi(f) F(\id,\id)g \\
\end{aligned}
\end{equation}
where the last step used the fact that $\sigma$ is a right $A$-module map. Combining \eqref{eq:BODGER} and \eqref{eq:BADGER} yields
\[   (\dif s_1 F+s_2\dif F)(f,g)
	 = F(f,g)  + \wtild{F}\sigma\left[ -\varphi(f)\ssp(g) + \ssp(fg) - \ssp(f)g \right] \;,  \]
so our proof is completed by observing that, for all $f,g\in A$\/,
\[ 
 \ssp(fg) 
 = fg -\varphi(f)\varphi(g)\id 
	   = \varphi(f)(g-\varphi(g)\id) + (f-\varphi(f)\id)g 
	   = \varphi(f)\ssp(g)+\ssp(f)g\;. 
\] 
\end{YCproof}

\begin{coroll}\label{c:proj-ideal-gives-splitting}
Let $A$\/, $\varphi$\/, $X$ and $C$ be as in Proposition \ref{p:proj-gen-split}. Then for any $n\geq 2$\/, there exist bounded linear maps $s_{n-1}:\Co{C}{n}(A,X)\to\Co{C}{n-1}(A,X)$ and 
 $s_n:\Co{C}{n+1}(A,X)\to\Co{C}{n}(A,X)$\/, such that $\norm{s_{n-1}}\leq 2C+1$\/, $\norm{s_n}\leq 2C+1$ and $s_n\dif+\dif s_{n-1}=\sid$\/.
 \end{coroll}

\begin{YCproof}
We use a `reduction-of-dimension' argument. By~\cite[\S1.a]{BEJ_CIBA}, for each fixed $n\geq 2$ there is a Banach $A$-bimodule $V_{n-2}$ (depending on $A$ and $X$), such that the cochain complex $\Co{C}{*}(A,V_{n-2})$ is \emph{isometrically} isomorphic to the cochain complex $\Co{C}{n-2+*}(A,X)$\/. Moreover, inspection of the relevant formulas in \cite[\S1.a]{BEJ_CIBA} shows that if the left action of $A$ on $X$ is given by $\varphi$\/, then so is the left action on $V_{n-2}$\/.
 Hence, by Proposition \ref{p:proj-gen-split} there exist $s_1:\Co{C}{2}(A,V_{n-2})\to\Co{C}{1}(A,V_{n-2})$ and $s_2:\Co{C}{3}(A,V_{n-2})\to\Co{C}{2}(A,V_{n-2})$ such that $\norm{s_1}\leq 2C+1$\/, $\norm{s_2}\leq 2C+1$ and $\dif s_1+s_2\dif =\sid$\/. Transferring these maps over to the complex $\Co{C}{n-2+*}(A,X)$, we obtain $s_{n-1}$ and $s_n$ as required. 
\end{YCproof}

\begin{rem}
Corollary \ref{c:proj-ideal-gives-splitting} does not guarantee a \emph{sequence} of maps $(s_n)_{n\geq 1}$ which satisfies $\dif s_n + s_{n+1}\dif =\sid$ for all $n$\/; it requires us to fix $n$ first. However, we could construct such a sequence, \emph{with the same constants as in Corollary~\ref{c:proj-ideal-gives-splitting}}, by taking
\[ (s_nT)(f_1,\ldots, f_n) = -\wtild{T}(\sigma\ssp(f_1)\tp f_2\tp\ldots f_n) + \varphi(f)T(\id,\id,f_2,\ldots, f_n) \;.
\]
Then, by a similar but more cumbersome calculation to that used in proving Proposition \ref{p:proj-gen-split}, it can be shown that $\dif s_n + s_{n+1}\dif=\sid$ for all $n$\/.
\end{rem}

We can weaken the hypotheses
of Proposition~\ref{p:proj-gen-split} to obtain a stronger technical result. This is not needed for the case of $\lp{1}(\Z_+)$ but will be used later for the other examples.

\begin{propn}\label{p:flat-gen-split}
Let $A$ be a commutative Banach algebra and let $\varphi\in\Phi_A$\/.
 Let $Y$ be a
Banach $A$-bimodule where the right action of $A$ is given by
\[  y\cdot f = \varphi(f)y \qquad(f\in A, y\in Y). \]
 Suppose that $A$ has a bounded approximate identity of bound $1$\/, and that $M_\varphi$ is $A$-flat with constant~$\leq C$\/.
Then there exist bounded linear maps $s_1:\Co{C}{2}(A,Y')\to \Co{C}{1}(A,Y')$ and $s_2:\Co{C}{3}(A,Y')\to\Co{C}{2}(A,Y')$ 
such that $s_2\dif+\dif s_1=\sid$\/.
\[ 
\begin{diagram}[tight,width=4em]
\Co{C}{1}(A,Y') & \RL{\dif}{s_1} & \Co{C}{2}(A,Y') & \RL{\dif}{s_2} & \Co{C}{2}(A,Y')
\end{diagram}
\]
Moreover, we can choose $s_1$ and $s_2$ such that $\max(\norm{s_1}, \norm{s_2}) \leq 2C+1$\/.
\end{propn}

We shall skip the proof as it is a more elaborate, and notationally more opaque, version of the proof of Proposition~\ref{p:proj-gen-split}. For sake of completeness, however, a full proof of Proposition~\ref{p:flat-gen-split} is given in an appendix.

As in the projective case, we can use a dimension-shift argument to extend this result, and obtain splitting maps in each degree $\geq 2$\/. The extension process is almost identical to that in the proof of Corollary~\ref{c:proj-ideal-gives-splitting}, so we omit the details.
\subsection*{The case of peak points}
Suppose $\Cplx_\varphi$ is $A$-flat: this is equivalent to requiring that $M_\varphi$ have a bounded approximate identity. In this case, the point cohomology groups $\Co{H}{n}(A,\Cplx_\varphi)$ are zero for all $n\geq 1$\/, and one can construct
splitting maps for the Hochschild chain complex in terms of the aforementioned bounded approximate identity in $M_\varphi$\/.

To our knowledge the formulas for these splitting maps are not spelled out explicitly in the literature, so for sake of completeness we give a quick account. (See, however, \cite[Proposition 4.19]{MCW_UA} and the remarks that follow for `implicit' statements.)

\begin{defn}
Let $\varphi$ be a character on $A$\/. A \dt{$\delta$-net for $\varphi$} is a net $(v_i)$ in $A$ with the property that $av_i\to 0$ and $v_ia\to 0$ for every $a\in M_\varphi$\/, and $\varphi(v_i)=1$ for all~$i$\/. If we furthermore have $\norm{v_i}\leq C$ for all $i$ then we say that our $\delta$-net is \dt{of bound~$C$}\/. 
\end{defn}

\begin{propn}\label{p:delta-net-gives-splitting}
Let $Y$ be a Banach $A$-bimodule (not necessarily contractive) where the right action of $A$ is given by
\[  y\cdot f = \varphi(f)y \qquad(f\in A, y\in Y). \]
Suppose that there exists a $\delta$-net for $\varphi$ of bound~$C$\/, and let $u\in A$ be such that $\varphi(u)=1$\/. Then there exist bounded linear maps $s_0:\Co{C}{1}(A,Y')\to \Co{C}{0}(A,Y')$ and $s_1:\Co{C}{2}(A,Y')\to\Co{C}{1}(A,Y')$\/, such that $s_1\dif+\dif s_0=\sid$\/.
\[ \begin{diagram}[tight,width=4em]
\Co{C}{0}(A,Y') & \RL{\dif}{s_1} & \Co{C}{1}(A,Y') & \RL{\dif}{s_1} & \Co{C}{2}(A,Y')
\end{diagram}\]
Moreover, we can choose $s_0$ and $s_1$ such that $\max(\norm{s_1}, \norm{s_2}) \leq C\norm{u}$\/.
\end{propn}

\begin{YCproof}
By hypothesis there exists a net $(v_\al)$ in $A$ such that $\norm{v_\al}\leq C$ and $\varphi(v_\al)=1$ for all~$\al$\/, and $av_\al, v_\al a\to 0$ for all $a\in M_\varphi$\/.
Our hypothesis on $Y$ ensures that:
\begin{itemize}
\item[(1)] for all $y\in Y$ and all $\al$\/, $y\cdot v_\al u = y$\/;
\item[(2)] for all $a\in A$\/, $v_\al(u\varphi(a)-ua)\to 0$\/.
\end{itemize}

Given $H\in \Co{C}{1}(A,Y')$ and $y\in Y$\/, put $(s_0^\al H)(y) = H(v_\al u)(y)$;
and given $F\in\Co{C}{2}(A,Y')$\/, $a\in A$ and $y\in Y$\/, put
$(s_1^\al F)(a)(y) = F(v_\al u, a)(y)$.
Then $s_0^\al$ and $s_1^\al$ are bounded linear maps with norm $\leq C\norm{u}$.
We have
\[  \dif s_0^\al H(a)(y)
	 = (s_0^\al H)(y\cdot a - a\cdot y)
	= (s_0^\al H) (\varphi(a)y-a\cdot y)  = H(v_\al u)(\varphi(a)y)-H(v_\al u)(a\cdot y)\,,\]
and, using observation (1) from earlier,
\[
 \begin{aligned}
s_1^\al\dif H(a)(y) = \dif H(v_\al u, a)(y) 
	& = H(a)(y\cdot v_\al u) - H(v_\al ua)(y) + H(v_\al u)(a\cdot y) \\
	& = H(a)(y) - H(v_\al ua)(y) + H(v_\al u)(a\cdot y) \;.
\end{aligned}
 \]
Hence
$(\dif s_0^\al H+s_1^\al\dif H)(a)(y)= H(a)(y) +  H(v_\al u\varphi(a)-v_\al ua)(y)$
for all $a\in A$\/, $y\in Y$\/.

Take $s_0$ to be a weak$^*$-cluster point of the net $(s_0^\al)$ in $\Lin{\Lin{A}{Y'}}{Y'}$, and $s_1$ to be a weak$^*$-cluster point of the net $(s_1^\al)$ in $\Lin{\Co{C}{2}(A,Y')}{\Co{C}{1}(A,Y')}$\/.
Using
our earlier observation~(2),
we find that $\dif s_0 + s_1\dif=\sid$ as required. The claimed upper bounds on $\norm{s_0}$ and $\norm{s_1}$ are immediate.
\end{YCproof}

As before, a dimension-shift argument extends the previous result to produce splitting maps in each degree. We omit the details.

%

\end{section}

\begin{section}{Application to $\lp{1}(\Z_+)$}\label{s:l1Z+}
In this section it is convenient to work with $\ATpos$\/: all the calculations that follow could instead be stated for $\lp{1}(\Z_+)$ by inverting Gelfand transforms, but the function-algebra perspective
is perhaps more intuitive.

If $\lm\in\clos{\Disc}\/,$ let $M_\lm=\{f\in \ATpos \st f(\lm)=0\}$\/.
For later reference we note that
\begin{equation}\label{eq:BluePeter}
\cG^{-1}(M_\lm)=\{ a\in \lp{1}(\Z_+) \st \sum_{n\in\Z_+} a(n) \lm^n = 0 \}\;.
\end{equation}

We first dispose of the case $\abs{\lm}=1$\/, i.e.~the case where our character is the restriction of a character on $\AT$\/. By Corollary~\ref{c:making-delta-nets} and taking Fourier transforms, there exists a $\delta$-net for this character which has bound~$1$.
Hence, by Proposition~\ref{p:delta-net-gives-splitting} and the remarks that follow~it, we obtain \emph{contractive} splitting maps in each degree.

We henceforth restrict attention to the case $\abs{\lm}<1$\/. As stated in the introduction, 
 the following result appears to be new.

\begin{thm}\label{t:main_result}
Let $\abs{\lm}<1$ and $n\geq 1$\/. Then there exist bounded linear maps $s_n$ and $s_{n+1}$ as shown in the diagram below,
such that $\norm{s_n}\leq 19$, $\norm{s_{n+1}}\leq 19$ and $\dif s_n+s_{n+1}\dif=\sid$\/.
\[ 
\begin{diagram}[tight,width=5em]
 \Co{C}{n}(\ATpos,\Cplx_\lm) & \RL{\dif}{s_n} & \Co{C}{n+1}(\ATpos,\Cplx_\lm) & \RL{\dif}{s_{n+1}} & \Co{C}{n+2}(\ATpos,\Cplx_\lm)
\end{diagram}
 \]
\end{thm}

As an immediate corollary, we get an affirmative answer to Question~\ref{mainQ}\/.
\begin{coroll}[Stability of point cocycles]
\label{c:stable-cocycle}
Let $n\geq 1$\/, let $\abs{\lm}<1$\/, and let $T\in\Co{C}{n}(\ATpos,\Cplx_\lm)$ be such that $\norm{\dif T}\leq 1$\/. Then there exists $S\in\Co{Z}{n}(\ATpos,\Cplx_\lm)$ such that $\norm{S-T}\leq 19$\/.
\end{coroll}

\begin{YCproof}
Let $s_n$\/, $s_{n+1}$ be as in Theorem \ref{t:main_result}. Then $\dif T = (\dif s_n +s_{n+1}\dif)\dif T = \dif s_n\dif T$\/. Put $S=T-s_n\dif T$\/.
\end{YCproof}

Theorem~\ref{t:main_result} is proved by combining Corollary~\ref{c:proj-ideal-gives-splitting}
with uniform bounds on the projectivity constant of $M_\lm$\/.
 The observation that $M_\lm$ is $\ATpos$-projective seems to be folklore, with the proof as follows: given $f\in M_\lm$\/, direct calculation with Fourier series shows that $(z-\lm)^{-1}f\in \ATpos$\/; hence the operator $\mathbf{s}:M_\lm \to M_\lm\ptp \ATpos$ that is given by $f\mapsto (z-\lm)\tp (z-\lm)^{-1}f$\/ satisfies the properties of Definition~\ref{dfn:proj-with-constant}.

By considering $\norm{\mathbf{s}(z^n-\lm^n)}$ for large $n$\/, it can be shown that $\norm{\mathbf{s}}\to\infty$ as $\abs{\lm}\nearrow 1$\/: thus $\mathbf{s}$ is insufficient for our purposes.
(We remark
that the splitting maps which would be produced by $\mathbf{s}$ on applying Proposition~\ref{p:proj-gen-split}, coincide with those implicitly obtained for the disc algebra by Johnson: see the proof of \cite[Proposition 9.1]{BEJ_CIBA}.)
Informally, our problems are caused by the fact that `dividing by $z-\lm$ increases the norm by a large factor'.
These problems disappear if
 we divide not by $z-\lm$, but by
the function
$b_\lm(z) = (z-\lm)(1-\lmbar z)^{-1}$\/.

\begin{lem} Let $\abs{\lm}<1$\/. Then: $b_\lm\in \ATpos$\/;
$b_\lm$ is invertible in $A(\Torus)$\/;
and
\[ \norm{b_\lm} = \norm{b_\lm^{-1}}_{\AT} =  1+ 2\abs{\lm}\;.\] 
\end{lem} 
\begin{YCproof}
For all $z\in\Disc$ we have
\[ b_\lm(z) = \lmbar^{-1} \left( \frac{1-\abs{\lm}^2}{1-\lmbar z} - 1 \right) = -\lm + \sum_{n\geq 1} (1-\abs{\lm}^2)\lmbar^{n-1} z^n \,,\]
so that $\norm{b_\lm}= \abs{\lm} + (1-\abs{\lm}^2)(1-\abs{\lm})^{-1}=1+2\abs{\lm}$\/. In particular, $b_\lm$ is an analytic function on the disc with absolutely convergent Taylor series.
The remaining observations follow immediately from the identity 
\[ b_\lm(e^{i\theta})^{-1} = \frac{1-\lmbar e^{i\theta}}{e^{i\theta}-\lm}
= \frac{e^{-i\theta}-\lmbar}{1-\lm e^{-i\theta}} = b_{\lmbar}(e^{-i\theta})\qquad(\theta\in\Real).\]
\smallskip
\end{YCproof}

Let $\ssS_\lm: A(\Torus)\to A(\Torus)$ be the bounded operator defined by
$\ssS_\lm(f) =  b_\lm^{-1}f$\/.
By the previous lemma, $\ssS_\lm$ is well-defined and invertible, with $\norm{\ssS_\lm}=\norm{\ssS_\lm^{-1}} = 1+2\abs{\lm}$\/. It is also clearly a (right) $A(\Torus)$-module map.

\begin{lem}\label{l:divisible}
$\ssS_\lm(M_\lm) \subseteq \ATpos$\/.
\end{lem}

\begin{YCproof}
Since $M_\lm$ is the closed linear span of the sequence $(z^n-\lm^n)_{n\geq 1}$\/, it suffices by linearity and continuity to show that $\ssS_\lm(z^n-\lm^n)\in \ATpos$ for all $n\in\Nat$\/. But this is obvious since
\begin{equation}\label{eq:reusable}
\ssS_\lm(e^{in\theta}-\lm^n) = \frac{(1-\lmbar e^{i\theta})(e^{in\theta}-\lm^n)}{e^{i\theta}-\lm} = (1-\lmbar e^{i\theta})\sum_{j=0}^{n-1} \lm^{n-1-j}e^{ij\theta}\;.
\end{equation}
\end{YCproof}


\begin{propn}\label{p:l1Z+_projmax}
If $\abs{\lm}<1$\/, then $M_\lm$ is $\lp{1}(\Z_+)$-projective with constant $\leq (1+2\abs{\lm})^2$\/.
\end{propn}

\begin{YCproof}
Define $\sigma_\lm: M_\lm\to M_\lm\ptp \ATpos$ by
$\sigma_\lm(f) = b_\lm\tp \ssS_\lm(f)$\/.
Clearly $\sigma_\lm$ is a bounded linear, right $\ATpos$-module map (since $\ssS_\lm$ is), and it satisfies
$\pi\sigma_\lm(f)= b_\lm b_\lm^{-1}f =f$
for all $f\in M_\lm$\/.
We finish by observing that $\norm{\sigma_\lm} \leq \norm{b_\lm}\norm{\ssS_\lm} \leq (1+2\abs{\lm})^2$\/.
\end{YCproof}

Theorem \ref{t:main_result} now follows from Proposition \ref{p:l1Z+_projmax} and Corollary~\ref{c:proj-ideal-gives-splitting}.

\begin{rem}
In view of Johnson's explicit formulas for splitting maps, we note that by specializing the proofs of the relevant results above, we obtain similar formulas.
Namely, in Theorem~\ref{t:main_result} we can take
\begin{equation}\label{eq:explicit-better-split}
 s_n T(f_1,\ldots, f_n)
 = \left\{ \begin{gathered}
    -  T( b_\lm, b_\lm^{-1}(f_1-f_1(\lm)\id), f_2, \ldots, f_n) \\
    +  f_1(\lm)T(\id,\id,f_2,\ldots, f_n)
\end{gathered}\right.
\end{equation}
for $f_1,\ldots, f_n\in\ATpos$\/.
Using the earlier calculation \eqref{eq:reusable}, it is straightforward if a little tedious to check that
\[ \begin{aligned}
\sup\{ \norm{b_\lm^{-1}(f-f(\lm)\id)} \st \norm{f}\leq 1\}
& = \sup_n \Norm{(\id -\lmbar z)(\sum\nolimits_{j=0}^{n-1}\lm^{n-1-j}z^j)} \\
& = 1+2\abs{\lm} 
\end{aligned}  \]
It follows that in Theorem~\ref{t:main_result} and Corollary~\ref{c:stable-cocycle}, the constant $19$ can be replaced throughout by~$10$\/.
The author does not know if this is optimal.

We note also that the formula \eqref{eq:explicit-better-split} also gives a well-defined, bounded linear map $\Co{C}{n+1}(A(\Disc),\Cplx_\lm) \to \Co{C}{n}(A(\Disc),\Cplx_\lm)$\/. Thus we obtain new splitting maps for the disc algebra case, different from those considered by Johnson. 
\end{rem}

\end{section}
\begin{section}{Further examples}
%
In this section we turn, as promised, to algebras of the form $\lp{1}(\bbG_+)$ where $\bbG\leq \Rat$\/, and to the algebra $L^1(\Real_+)$. For sake of brevity, we
 shall only verify that the flatness constants for maximal ideals in these algebras are uniformly bounded; just as for $\lp{1}(\Z_+)$\/, one can then go on to deduce splitting results with uniform control for the corresponding point cohomology, using Proposition~\ref{p:flat-gen-split} instead of Proposition~\ref{p:proj-gen-split}.

\begin{subsection}{The case of $\ell^1(\bbG_+)$ for $\bbG\leq \Rat$}
\label{ss:l1Q+}
Fix a subgroup $\bbG\leq \Rat$
such that $\bbG$ is dense in $\Real$\/. (It is an easy exercise to show that any non-dense subgroup of $\Real$ is of the form $\al\Z$ for some $\al>0$\/). By the results of~\cite[\S4]{AreSin_gaf56}, the characters of $\ell^1(\bbG_+)$ have one of two forms:
\begin{enumerate}
\item $b \mapsto \sum_{x\in\bbG_+} b(x) e^{-tx}\chi(x)$\/, 
 where $\chi\in\dual{\bbG}$ and $0\leq t <\infty$\/;
\item $b\mapsto b(0)$\/.
\end{enumerate}
Informally, the second case corresponds to taking `$t=+\infty$'\/.

Any character for which $t=0$ extends uniquely to a character of the group algebra $\lp{1}(\bbG)$\/. By Corollary~\ref{c:making-delta-nets}, for each such character there exists a $\delta$-net of bound~$1$\/; hence
by Proposition~\ref{p:delta-net-gives-splitting} and the remarks that follow~it, we obtain \emph{contractive} splitting maps in each degree.

We now restrict attention to the case $0< t <\infty$\/; the `$t=+\infty$' case will be treated separately, albeit with sketches rather than full proofs. Thus, fix $t\in(0,\infty)$ and $\chi\in\dual{\bbG}$\/, and let $\varphi$ be the character on $\lp{1}(\bbG_+)$ defined by
$\varphi(b) = \sum_{x\in \bbG_+} b(x) e^{-tx} \chi(x)$\/.
Let $\bM_\varphi=\ker\varphi$\/. We shall construct a net $(\rho_\al)_{\al\in\cD}$ of bounded linear maps $\bM_\varphi\to \bM_\varphi\ptp\lp{1}(\bbG_+)$, which satisfies the conditions of Lemma \ref{l:flat-by-approx} and is such that $\sup_\al \norm{\rho_\al} \leq 9$\/. The indexing set $\cD$ is defined as follows: let ${\bbG}_{>0}$ denote the set $\{ \al\in\bbG \st \al >0 \}$\/, and equip $\bbG_{>0}$ with the partial order given by
\[ \text{$\al\preceq\beta$ if and only if $\al=k\beta$ for some positive integer~$k$\/;}\]
now define $\cD$ to be the partially ordered set $( \bbG_{>0},\preceq)$\/. To use $\cD$ to index our net, we need a small but important lemma.

\begin{lem}\label{l:HCF}
$\cD$ is a filtered set.
\end{lem}
\begin{YCproof}
Let $\al_1$\/, $\al_2\in\cD$\/.
%
Since $\al_1$ and $\al_2$ are positive rationals, by clearing denominators there exist $\gm\in\Rat_{>0}$ and coprime, positive integers $n_1$ and $n_2$\/, such that $\al_i=n_i\gm$\/, $i=1,2$\/.
Then by the Euclidean algorithm there exist integers $b_1, b_2$ such that $b_1n_1+b_2n_2=1$\/. Hence $\gm= b_1n_1\gm + b_2n_2\gm = b_1\al_1 + b_2\al_2 \in \bbG$\/, and so $\gm$ is a common upper bound in $\cD$ for $\al_1$ and $\al_2$\/, as required.
\end{YCproof}

For each $\al\in\cD$ define $\theta_\al:\lp{1}(\bbG_+)\to \lp{1}(\bbG_+)$ by
\[ \theta_\al(b) = \sum_{n\in\Z_+} \left( \sum_{x\in \bbG, n\al\leq x < (n+1)\al}   b(x)e^{-t(x-n\al)}\chi(x-n\al) \right) \delta_{n\al} \;.\] 
Then $\theta_\al$ is a \emph{contractive} linear \emph{projection} of $\lp{1}(\bbG_+)$ onto the closed subalgebra $\lp{1}(\al\Z_+)$\/.

\begin{lem}\label{l:discretize}
For every $b\in\lp{1}(\bbG_+)$ we have $b=\lim_{\al\in\cD} \theta_\al(b)$\/.
\end{lem}

\begin{YCproof}
Let $\veps>0$\/. There exists a finite set $F\subset \bbG_+$ such that
$\sum_{x\in \bbG_+\setminus F} \abs{ b(x) } < \veps/2$\,: set $b_0\defeq \sum_{x\in F} b(x)\delta_x$\/.
Since $F$ is finite and $\cD$ is filtered, there exists an upper bound $\al_0\in\cD$ for all elements of $F\setminus\{0\}$. Then, for all $\al\succeq\al_0$ in~$\cD$\/, we have $F\subset\al\Z_+$ and so $\theta_\al(b_0)=b_0$\/.
Hence for all $\al\succeq\al_0$\/, $\norm{\theta_\al(b)-b} \leq \norm{\theta_\al(b-b_0)-(b-b_0)} \leq\veps$\,.
\end{YCproof}

\begin{lem}\label{l:staying-in-ideal}
If $b\in \bM_\varphi$ then $\theta_\al(b)\in \bM_\varphi\cap\lp{1}(\al\Z_+)$\/.
\end{lem}

\begin{YCproof}
\[ \begin{aligned}
   \varphi(\theta_\al(b)) 
  & = \sum_{y\in \bbG_+} (\theta_\al b)(y)e^{-ty}\chi(y) \\
  & = 
\sum_{n\in\Z_+} (\theta_\al b)(n\al)e^{-tn\al}\chi(n\al) \\
  & =
\sum_{n\in\Z_+} \left( \sum_{x\in \bbG, n\al\leq x < (n+1)\al} b(x)e^{-t(x-n\al)}\chi(x-n\al)\right) e^{-tn\al}\chi(n\al) \\
  & = \sum_{x\in \bbG_+} b(x)e^{-tx}\chi(x) = 0 \;.
\end{aligned} \]
\end{YCproof}

\begin{propn}\label{p:Rat-flat-ideal}
$\bM_\varphi$ is $\lp{1}(\bbG_+)$-flat with constant~$\leq 9$\/.
\end{propn}

\begin{YCproof}
Let $\al\in\cD$ and put $\bM^{(\al)}_\varphi=\bM_\varphi\cap\lp{1}(\al\Z_+)$\/.
Then
\[ 
\bM^{(\al)}_\varphi
	 = \{ b \in \lp{1}(\al\Z_+) \st \varphi(b)=0 \} 
	 = \left\{ b\in\lp{1}(\al\Z_+) \st \sum_{n\in \Z_+} b(n\al) \bigl(e^{-t\al}\chi(\al)\bigr)^n = 0 \right\}\;.
\]
Comparing this with \eqref{eq:BluePeter} above, 
we see that $\bM^{(\al)}_\varphi$ is isometrically isomorphic as a Banach space to $M_\lm$ where $\lm=e^{-t\al}\chi(\al)$\/, and that this isomorphism intertwines the action of $\lp{1}(\al\Z_+)$ on $\bM^{(\al)}_\varphi$ with the action of $A(\Torus)^+$ on $M_\lm$\/. Hence by Proposition~\ref{p:l1Z+_projmax}, there exists a bounded linear, right $\lp{1}(\al\Z_+)$-module map $\sigma_\al: \bM^{(\al)}_\varphi \to \bM^{(\al)}_\varphi\ptp\lp{1}(\al\Z_+)$ such that
\[ \pi\sigma_\al(f)=f \qquad\text{ for all $f\in \bM^{(\al)}_\varphi$\/}  \]
and such that $\norm{\sigma_\al} \leq (1+2 \abs{e^{-t\al}\chi(\al)})^2 \leq 9$\/.
Via the inclusion of $\bM^{(\al)}_\varphi\ptp\lp{1}(\al\Z_+)$ into $\bM_\varphi\ptp\lp{1}(\bbG_+)$\/, we can regard $\sigma_\al$ as taking values in $\bM_\varphi\ptp\lp{1}(\bbG_+)$\/. By Lemma~\ref{l:staying-in-ideal}, $\theta_\al(\bM_\varphi)\subseteq \bM^{(\al)}_\varphi$ and we may therefore set
\[ \rho_\al=\sigma_\al\theta_\al : \bM_\varphi \to \bM_\varphi \ptp\lp{1}(\bbG_+) \;. \]

It suffices to show that the net $(\rho_\al)$ satisfies conditions $(i)$ and $(ii)$ in Lemma~\ref{l:flat-by-approx}.
Firstly, for every $c\in \bM_\varphi$\/,
$\lim_\al \pi\rho_\al(c) = \lim_\al \pi\sigma_\al(\theta_\al c) =  \lim_\al\theta_\al(c) =c$\/;
this shows that $(i)$ holds. 
Secondly, let $c\in \bM_\varphi$, $b\in\lp{1}(\bbG_+)$\/.
Then,
for each $\al\in\cD$\/, since $\sigma_\al$ is a right $\lp{1}(\al\Z_+)$-module map, we have
$\sigma_\al\theta_\al(c)\cdot\theta_\al(b) = \sigma_\al\bigl( \theta_\al(c)\theta_\al(b) \bigr)$.
Hence
\[ \begin{aligned}
  \norm{\rho_\al(c)\cdot b - \rho_\al(cb) }
    & = \norm{ \sigma_\al\theta_\al(c)\cdot b - \sigma_\al\theta_\al(cb) } \\
    & = \norm{ \sigma_\al\theta_\al(c)\cdot(b-\theta_\al(b))
		+ \sigma_\al\bigl( \theta_\al(c)\theta_\al(b) - \theta_\al(cb) \bigr) } \\
    & \leq KC\norm{c}\norm{b-\theta_\al(b)} + K\norm{\theta_\al(c)\theta_\al(b)-\theta_\al(cb) }
\end{aligned} \]
and by Lemma \ref{l:discretize}, the right hand side converges to $0$ as $\al$ `tends to infinity', so that $(ii)$~holds.
\end{YCproof}

This completes our analysis for those characters with $0<t<\infty$\/. It remains to deal with the maximal ideal
$\bM_{(\infty)} = \{ a\in\lp{1}(\bbG_+) \st a(0)=0 \}$ that
corr\-es\-ponds to the case `$t=+\infty$'. 
Most of the preceding analysis works in this case also, provided that one reinterprets certain formulas appropriately. For example, the approximating maps $\theta_\al: \lp{1}(\bbG_+)\to\lp{1}(\bbG_+)$ are now defined by
$\theta_\al(b) =\sum_{n\in\Z_+} b(n\al)\delta_{n\al}$\/.
In this particular case we can write down an explicit bounded net $(\rho_\al)_{\al\in\cD}$\/, given by
\[ \rho_\al(f) = \delta_\al \tp (\delta_{-\al}*\theta_\al(f)) \qquad(f\in \bM_{(\infty)}). \]
Since $\theta_\al(f)\to f$ for each $f\in\lp{1}(\bbG_+)$\/, it is easily checked that the net $(\rho_\al)_{\al\in\cD}$ satisfies the conditions of Lemma~\ref{l:flat-by-approx}, and we conclude that $\bM_{(\infty)}$ is $\lp{1}(\bbG_+)$-flat with constant~$1$\/.
\end{subsection}

\begin{subsection}{The case of $L^1(\Real_+)$}\label{ss:L1R+}
The results of this section could be obtained by approximation arguments similar to those used in passing from $\lp{1}(\Z_+)$ to $\lp{1}(\bbG_+)$\/. However,
 it is simpler to write down continuous analogues of the constructions in Section~\ref{s:l1Z+}, aided by a little guesswork and hindsight.

Recall from Section~\ref{s:L1prelim} that the character space of $L^1(\Real_+)$ may be identified with~$\clos{\UHP}$\/, and that those characters corresponding to the boundary of $\clos{\UHP}$ are restrictions of characters on $L^1(\Real)$\/. As with our previous examples, Corollary~\ref{c:making-delta-nets} ensures that for each such character there exists a $\delta$-net in $L^1(\Real_+)$ of bound~$1$\/. Hence, by Proposition~\ref{p:delta-net-gives-splitting} and the remarks that follow~it, we obtain \emph{contractive} splitting maps in each degree.

For the rest of this section, we therefore restrict ourselves to those point modules which correspond to points of the \emph{open} upper-half plane~$\UHP$\/. Fix $\lm \in\UHP$\/, and define $\upos{\lm}$ and $\uneg{\lm}$ in $L^1(\Real)$
 by
$\upos{\lm}(t) = \indic{t > 0} e^{-i\lmbar t}$\/,
\;$\uneg{\lm}(t) = \indic{t < 0} e^{-i\lm t}$\/.
Since $\Im\lm > 0$\/, $\upos{\lm}$ and $\uneg{\lm}$ both lie in
$L^1(\Real)$\/. Moreover
\begin{equation}\label{eq:a-norm-bound}
\norm{\upos{\lm}} = \norm{\uneg{\lm}} = \int_0^\infty e^{-(\Im \lm )t}\,dt = \frac{1}{\Im\lm} \;.
\end{equation}
The Fourier transforms of both $\upos{\lm}$ and $\uneg{\lm}$ are easily computed: they are
\begin{equation}\label{eq:FTs}
 \cF\upos{\lm}(x) = \frac{i}{x-\lmbar} \quad,\quad
\cF \uneg{\lm}(x) = \frac{1}{i(x-\lm)} \quad(x\in\Real).
\end{equation}


Set $h_\lm \defeq \delta_0 - 2 (\Im\lm) \upos{\lm}$ and $\check{h}_\lm \defeq \delta_0 - 2(\Im\lm )\uneg{\lm}$\/, regarding both as elements of the measure algebra $M(\Real)$\/. By \eqref{eq:a-norm-bound} we have $\norm{h_\lm}=\norm{\check{h}_\lm}=3$\/, and by \eqref{eq:FTs} we have
\begin{equation}\label{eq:BEAKER}
 \cF h_\lm(x) = 1 -\frac{2i\Im\lm}{x-\lmbar} = \frac{x-\lm}{x-\lmbar} \quad\text{and}\quad \cF \check{h}(x) = 1 + \frac{2i\Im\lm}{x-\lm} = \frac{x-\lmbar}{x-\lm} \;.
\end{equation}

\begin{coroll}\label{c:inverse}
Let $f\in L^1(\Real)$\/. Then $h_\lm*\check{h}_\lm*f=f$ \textup{(}a.e.\textup{)}
\end{coroll}

\begin{YCproof}
One could check this by direct computation with convolutions. More intuitive is the following argument: by \eqref{eq:BEAKER} we have
\[ \cF( h_\lm*\check{h}_\lm*f-f)= \cF h_\lm \cF\check{h}_\lm\cF f-\cF f=0\;; \]
and the corollary now follows from the injectivity of
$\cF: L^1(\Real)\to C_0(\Real)$ (see \cite{DJN_unique-ft} for a particularly slick proof).
%
\end{YCproof}

The function $\cF h_\lm$ plays an analogous role to the function $b_\lm$ from Section~\ref{s:l1Z+}, and the following lemma may be seen as a `continuous analogue' of Lemma~\ref{l:divisible}.
\begin{lem}
Let $f\in \cM_\lm$\/. Then $(\check{h}_\lm*f)(t)= 0$ for a.e.~$t<0$\/, so that we may identify $\check{h}_\lm*f$ with a (unique) element of $L^1(\Real_+)$\/.
\end{lem}

\begin{YCproof}
For $t<0$ we have
\[ \begin{aligned}
(\check{h}_\lm*f)(t) = -2(\Im\lm) (\uneg{\lm}*f)(t)
	& = -2(\Im\lm) \int_{-\infty}^0 e^{-i\lm s} f(t-s)\,ds \\
	& = 2(\Im\lm) \int_t^\infty e^{i\lm (r-t)} f(r)\,dr \\
	& = 2(\Im\lm) e^{-i\lm t} \left[ \int_t^0 e^{i\lm r} f(r)\, dr + \cF f(\lm) \right]
	& = 0\;,
 \end{aligned} \]
the last step following since $f=0$ a.e.~on the interval $[t,0]$ and $\cF f(\lm)=0$\/. 
\end{YCproof}

\begin{coroll}
$\cM_\lm$ is $L^1(\Real_+)$-flat with constant less than or equal to~$9$\/.
\end{coroll}

\begin{YCproof}
Fix a bounded approximate identity $(e_\alpha)$ in $L^1(\Real_+)$ which satisfies $\norm{e_\alpha}=1$ for all~$\alpha$
(for instance, take normalized characteristic functions of small intervals containing~$0$).
Define $\rho_\al: \cM_\lm \to \cM_\lm\ptp L^1(\Real_+) \subset \cM_\lm\ptp\fu{L^1(\Real_+)}$ by
\[ \rho_\al(f) = (e_\al* h_\lm)\tp (\check{h}_\lm*f) \qquad(f\in\cM_\lm). \]
Then each $\rho_\al$ is a right $L^1(\Real_+)$-module map, satisfying $\norm{\rho_\al}\leq\norm{h_\lm}\norm{\check{h}_\lm}\leq 9$. By Corollary~\ref{c:inverse},
$\pi\rho_\al(f) =e_\al*h_\lm*\check{h}_\lm*f = e_\al*f$
for each $f\in\cM_\lm$\/. The result now follows by applying Lemma~\ref{l:flat-by-approx}.
\end{YCproof}

\end{subsection}
\end{section}

\begin{section}{Closing remarks}
We have seen, in particular, that those maximal ideals in $\lp{1}(\bbG_+)$ or $L^1(\Real_+)$ which are not in the (Shilov) boundary are flat. The author does not know if any of these ideals are projective, and if so whether there are uniform bounds on the projectivity constants.

The case $G=\Real_d$ is not covered by the results above: this is because Lemma~\ref{l:HCF} fails, so that we cannot approximate $\lp{1}((\Real_d)_+)$ by $\lp{1}(\Z_+)$ in such a way as to deduce flatness of maximal ideals. More abstractly, subgroups of $\Rat$ are rank one abelian groups, while $\Real_d$ is not. Lemma~\ref{l:HCF} also fails if we take $G=\Z+\Z\theta\subset\Real$ or $\Rat+ \Rat\theta\subset\Real$ for some fixed irrational, positive $\theta$\/; it seems that to treat these examples, additional approximation techniques will be needed.

A natural first step towards understanding higher rank cases is to consider the algebras $\lp{1}(\Z_+^k)$ where $2\leq k <\infty $\/. Here, the point modules corresponding to interior points of the maximal ideal space have non-trivial cohomology in degrees $k$ and below, and vanishing cohomology in degrees $k+1$ and above. Briefly: whereas in the one-variable case we had a short exact sequence
\[ 0 \leftarrow \Cplx_\varphi \leftarrow \lp{1}(\Z_+) \leftarrow M_\varphi \leftarrow 0 \]
in $k$ variables we have a \dt{projective resolution of length $k$}:
\[ 0 \leftarrow \Cplx_\varphi \leftarrow \lp{1}(\Z_+^k) \leftarrow P_1 \leftarrow \ldots P_k \leftarrow  0 \]
There is a canonical construction for the complex $P_*$ (namely, the \dt{Koszul resolution} associated to the module $\Cplx_\varphi$). If we follow this construction then each $P_i$ is a free $\lp{1}(\Z_+^k)$ module and hence is projective with constant~$1$\/. However, the norms of splitting maps for the corresponding Hochschild chain complex $\Co{C}{*}(\lp{1}(\Z_+^k),\Cplx_\varphi)$ will depend on the norms of the linear splitting maps for the resolution $P_*$\/, and our analysis therefore becomes more complicated. We hope to address this in future work.

\begin{acknowledgements}
The work presented here was motivated by discussions at the University of Newcastle-upon-Tyne, and carried out while the author was a postdoctoral researcher at the University of Manitoba. He thanks T.~G. Kucera for useful exchanges concerning the rank of abelian groups.

The article was prepared using Paul Taylor's \texttt{diagrams.sty} macros.
\end{acknowledgements}
\end{section}

\appendix

\begin{section}{Full proofs for the case of a flat maximal ideal}
To deal with the case where our maximal ideal $M_\varphi$ is merely flat, rather than projective, it is useful to set up some preliminary results on \wstar-continuous extensions.

\begin{lem}
\label{l:wstar-extns}
Let $E$, $F$ and $N$ be Banach spaces, and let $T\in\Lin{E\ptp F}{N'}$.
\begin{enumerate}
\item $T$ extends uniquely to a bounded linear, \wstar-\wstar-continuous map $\bT:(E\ptp F)''\to N'$\/, and $\norm{\bT}=\norm{T}$\/.
\item $T$ extends uniquely to a bounded linear map $\bT_{(1)}: E'' \ptp F \to N'$ that is \wstar-\wstar-continuous in the first variable, and $\norm{\bT_{(1)}}=\norm{T}$\/.
\item There exists a norm-one linear map ${\boldsymbol\kappa}_{(1)}: E''\ptp F\to (E\ptp F)''$ such that ${\boldsymbol\kappa}_{(1)}\circ(\kappa_E\tp\id) = \kappa_{(E\ptp F)}$\/.
\end{enumerate}
\end{lem}

\begin{YCproof}
For each $x\in N$\/, let $T(\blank)(x)$ denote the functional on $E\ptp F$ that is defined  by $e\tp f\mapsto T(e\tp f)(x)$\/; and for each $f\in F$ let $T(\blank\tp f)(x)$ denote the functional on $E$ that is defined by $e\mapsto T(e,f)(x)$\/. We then define $\bT$ and $\bT_{(1)}$ by
\begin{subequations}
\begin{flalign}
 \quad \bT({\mathbf h})(x) &= {\mathbf h}[T(\blank)(x)] &\text{for ${\mathbf h}\in (E\ptp F)''$ and $x\in N$\/,} \quad\\
 \quad \bT_{(1)}({\mathbf e}\tp f)(x) & = {\mathbf e}[T(\blank\tp f)(x)] & \text{for ${\mathbf e}\in E''$, $f\in F$ and $x\in N$\/.}\quad
\end{flalign}
\end{subequations}
It is easily checked that $\bT$ and $\bT_{(1)}$ have the required linearity, continuity and extension properties: uniqueness follows since a Banach space is \wstar-dense in its second dual. Finally, the last assertion follows by taking $N$ to be $(E\ptp F)''$ and $T$ to be $\kappa_{(E\ptp F)}$\/.
\end{YCproof}

\begin{lem}\label{l:flat-tweaked}
Let $A$ be a Banach algebra which has a bounded left approximate identity $(e_i)$ of bound~$K\geq 1$\/, and let $M$ be a right Banach $A$-module satisfying $\clos{MA}=M$\/. If $M$ is $A$-flat with constant $\leq C$\/, then there exists a bounded linear\, left $A$-module map $\Lm: (M\ptp A)'\to M'$ such that $\norm{\Lm}\leq CK$ and $\Lm\pi'(\psi)=\psi$ for all $\psi\in M'$\/.
\end{lem}

The lemma is a quantitative variant on a standard theme: we give the proof below, since we wish to avoid using the equivalence of our definition of flatness with the original one, and so cannot easily appeal to the arguments of \cite{Hel_HBTA} or~\cite{MCW_UA}.

\begin{YCproof}
By hypothesis (see Definition~\ref{dfn:flat} above) there exists a bounded linear\, left $A$-module map $\rho: (M\ptp\cu{A})'\to M'$ such that $\norm{\rho}\leq C$ and $\rho\pi'(\psi)=\psi$ for all $\psi\in M'$\/.
For each $i$\/, define a bounded linear map $\bj_i: (M\ptp A)'\to (M\ptp\cu{A})'$ by
\[ \bj_i \Psi (x\ptp c) = \Psi(x\ptp e_ic) \qquad(x\in M, c\in \cu{A}).\]
The net $(\bj_i)$ is bounded: let $\bj$ be a cluster point of this net in the \wstar-topology of ${\mathcal L}((M\ptp\cu{A})')$\/. Clearly $\norm{\bj}$ is a left $A$-module map, with norm $\leq K$.
Note that for each $\Psi\in (M\ptp A)'$ and $a\in A$\/,
\begin{equation}\label{eq:needed-shortly}
\bj(a\cdot\Psi)=a\cdot\Psi \,. \tag{$*$}
\end{equation}
(By linearity and continuity,  it suffices to check that both sides of \eqref{eq:needed-shortly} agree on elementary tensors of the form $x\tp c$\/, where $x\in M$ and $c\in \cu{A}$\/.
This in turn follows from the definition of the net $(\bj_i)$ and the fact that
 $\lim_i e_i(ca) = ca$\/.)

Put $\Lm=\rho\bj : (M\ptp A)' \to M'$\,: then $\Lm$ has norm $\leq CK$ and 
 is a left $A$-module map.
It remains only to show that $\Lm\pi'(\psi)=\psi$ for all $\psi\in M'$\/: since $MA$ is dense in $M$\/, it suffices to show that $\Lm\pi'(\psi)(xa)=\psi(xa)$ for all $\psi\in M'$\/, $x\in M$ and all $a\in A$\/. This follows from the following calculation:
\begin{flalign*}
\qquad & \Lm\pi'(\psi)(xa) - \psi(xa) \\ 
	& = \rho\bj\pi'(\psi)(xa) - \rho\pi'\psi(xa) & \text{(definitions of $\rho$ and $\Lm$)} \\
	& = \left(a\cdot\left[\rho\bj\pi'(\psi) - \rho\pi'(\psi)\right]\right)(x) & \text{(left action of $A$)} \\
	& = \rho\left[ \bj (a\cdot\pi'(\psi))- a\cdot\pi'(\psi)\right](x) & \text{($\rho$\/, $\bj$ are left $A$-module maps)} \\
	& = 0 &	\text{(by \eqref{eq:needed-shortly} above).}
\end{flalign*}
\end{YCproof}

\begin{YCnameproof}{Proof of Proposition~\ref{p:flat-gen-split}}
By rescaling the original bounded approximate identity, we may assume it has the form $(e_j)_{j\in\cJ}$ where $\varphi(e_j)=1$ for all $j$ and $\limsup_{j\in\cJ}\norm{e_j}=1$\/. Define $\ssp_i : A\to M_\varphi$ by
$\ssp_i(f) = f- \varphi(f)e_i$\/, so that for each $i$ $\ssp_i$ is a projection of norm $\leq 2$ of $A$ onto $M_\varphi$\/.

By Lemma~\ref{l:flat-tweaked} there exists a bounded linear, left $A$-module map $\Lm:(M_\varphi\ptp A)' \to M_\varphi'$ with $\norm{\Lm}\leq C$ and $\Lm\pi'(\psi)=\psi$ for all $\psi\in M_\varphi'$\/.
For each $i,j\in\cJ$\/, by using Lemma~\ref{l:wstar-extns}, we define $t^i_1, u^{i,j}_1:\Co{C}{2}(A,Y')\to\Co{C}{1}(A,Y')$ and
 $t^j_1, u^{i,j}_2:\Co{C}{3}(A,Y')\to\Co{C}{2}(A,Y')$ by
\[
\begin{aligned}
t^j_1 F(f) = - \wtild{\bF}\Lm'(\ssp_j f)
 \quad&,&\quad
t^j_2 G(f,g)  = - \wtild{\bG}_1(\Lm'(\ssp_j f)\tp g)
\\
u^{i,j}_1F(f)  = \varphi(f) F(e_i,e_j)
\quad &,&\quad
u^{i,j}_2G(f,g)  = \varphi(f) G(e_i,e_j,g)\;.
\end{aligned}\]
(In the definition of $t^j_2$\/, we are applying Lemma \ref{l:wstar-extns} with $E=M_\varphi\ptp A$\/, $F=A$ and $N=Y$\/.)

Recall the identity \eqref{eq:GOOFY} that was used earlier, which says that
\[ \wtild{\dif F}(w\tp b) = - \wtild{F}(\pi(w)\tp b) + \wtild{F}(w\cdot b) -\wtild{F}(w)b
\quad\text{for all $w\in M_\varphi\ptp A$\/ and $b\in A$\/.}\]
We may extend it by $w^*$-continuity, using the uniqueness parts of Lemma~\ref{l:wstar-extns}, to show that
{for all $b\in A$ and  ${\mathbf w}\in (M_\varphi\ptp A)''$\/,}
\begin{equation}\label{eq:beefed-up-GOOFY} 
\wtild{{\boldsymbol\dif}\bF}_{(1)}( {\mathbf w}\tp b) =
 - \wtild{\bF}_{(1)}(\pi''({\mathbf w})\tp b) + \wtild{\bF}({\mathbf w}\cdot b) - \wtild{\bF}({\mathbf w})b \,.
\end{equation}
Now if $b\in A$ and $h\in M_\varphi$ we have
$\wtild{\bF}_{(1)}(\pi''\Lm'(h)\tp b)=\wtild{\bF}_{(1)}(\kappa_A(h)\tp b) = F(h,b)$\/, since $\wtild{\bF}_1$ extends $\wtild{F}$\/. Combining this observation with~\eqref{eq:beefed-up-GOOFY} and the fact that $\Lm'$ is a right $A$-module map,
\[ \begin{aligned}
 t^j_2(\dif F)(f,g)
   & = \wtild{{\boldsymbol\dif}\bF}_{(1)}(\Lm'(\ssp_j f)\tp g) \\
   &  = \wtild{\bF}_{(1)}(\pi''\Lm'(\ssp_j f)\tp g) - \wtild{\bF}(\Lm'(\ssp_j f)\cdot g)  + (\wtild{\bF}\Lm'(\ssp_j f))g \\
   & =  F(\ssp_j f, g) - \wtild{\bF}\Lm'((\ssp_j f)g) + (\wtild{\bF}\Lm'(\ssp_j f))g\,.
\end{aligned} \]
Since we also have
\[ \begin{aligned}
   \dif(t^j_1 F)(f,g)
   &  = \varphi(f)t^j_1 F(g) - t^j_1 F(fg) + t^j_1 F(f)g \\
   & =   - \varphi(f)\wtild{\bF}\Lm'(\ssp_jg)
	  + \wtild{\bF}(\Lm'\ssp_j(fg)) - (\wtild{\bF}\Lm'(\ssp_j f))g \\
   & =  \wtild{\bF}\Lm'(fg-\varphi(f)g) - (\wtild{\bF}\Lm'(\ssp_jf))g \,,
\end{aligned} \]
we arrive at
\begin{equation}\label{eq:almost-split1}
\begin{aligned}
    (\dif t^j_1 F+t^j_2\dif) F(f,g)
	& = F(\ssp_j f, g)+\wtild{\bF}\Lm'(\varphi(f)(e_jg-g)) \\
	& = F(\ssp_j f, g) 
	+  \wtild{\bF}\Lm'( \varphi(f)(e_jg-g))  \\
	& = F(f, g) - \varphi(f)F(e_j,g) 
	+  \wtild{\bF}\Lm'(\varphi(f)(e_jg-g)) \;.
\end{aligned}
\end{equation}
On the other hand, a direct calculation (using the fact that $\varphi$ is multiplicative) yields
\[ \begin{aligned}
  \dif(u^{i,j}_1 F)(f,g)
	 & = \varphi(f) u^{i,j}_1F(g) - u^{i,j}_1F(fg)
		+ u^{i,j}_1F(f)g \\
	 & = \varphi(f) F(e_i,e_j)g\;,\quad\text{and}\\
    u^{i,j}_2(\dif F)(f,g) 
	& = \varphi(f) \dif F(e_i,e_j,g) \\
	& =
	\varphi(f) F(e_j,g)
	  	- \varphi(f)\left[ F(e_ie_j, g)- F(e_i, e_jg)\right] 
	 	 - \varphi(f) F(e_i,e_j)g\;,
\end{aligned} \]
so that
\begin{equation}\label{eq:almost-split2}
 (\dif u^{i,j}_1 + u^{i,j}_2\dif) F(f,g) = \varphi(f)F(e_j,g) -\varphi(f)\left[ F(e_ie_j,g) - F(e_i,e_jg)\right] \;.
\end{equation}

Putting $s^{i,j}_k=t^j_k+u^{i,j}_k$ for $k=1,2$\/, we deduce from \eqref{eq:almost-split1} and \eqref{eq:almost-split2} that
\[ (\dif s^{i,j}_1 +s^{i,j}_2\dif)F(f,g)
	 = F(f,g) + \varphi(f)\left[
		 \wtild{\bF}\Lm'(e_jg-g) + F(e_ie_j,g) -F(e_i, e_jg)
		\right] \;,
 \]
where for fixed $i\in\cJ$ and $f,g\in A$\/, the term in square brackets converges to zero as we take \wstar-limits with respect to~$j$\/. 
Also, $\norm{s^{i,j}_k}\leq 2C+\norm{e_i}\norm{e_j}$ for $k=1,2$\/.
Thus on taking
$s_k=\wstar\!\lim_i \wstar\!\lim_j s^{i,j}_k$
for $k=1,2$\/,
we obtain maps with the required properties, and the proof is complete.
\end{YCnameproof}
\end{section}



\vskip2.0em

\noindent{Y. Choi}

\noindent D\'epartement de math\'ematiques et de statistique, Universit\'e Laval, Qu\'ebec, Canada G\textup{1}V \textup{0}A\textup{6}

\noindent Email:\textup{: \texttt{y.choi.97@cantab.net}}

\smallskip\hfill Qu\'ebec City, October 2008

\hfill Revised May 2009
\end{document}